\documentclass[12pt]{article}

\usepackage{latexsym,amsmath, natbib, amsbsy, amsfonts, amssymb}

\newtheorem{theorem}{Theorem}
\newtheorem{proposition}{Proposition}
\newtheorem{remark}{Remark}
\def\text#1{\hbox{#1}}

\def\xib{\mbox{\boldmath$\xi$}}

\def\Rbb{{\mathbb R}}

\def\C{{\bf C}}
\def\Hb{{\bf H}}
\def\P{{\bf P}}
\def\Pb{{\bf P}}
\def\n{\noindent}

\def\iy{\infty}
\def\l{\lambda}

\def\xb{{\bf x}}

\def\x{{\bf x}}
\def\Xb{{\bf X}}

\def\s{\sigma}
\def\t{\theta}

\begin{document}

\title{An extremal problem with applications to testing multivariate independence}

\author{Alexander Nazarov\footnote{\em{Department of Mathematics and Mechanics, St.Petersburg State University, St. Petersburg, 198504, Russia}}
\ and \ Natalia Stepanova\footnote{\em{School of Mathematics and Statistics,
Carleton University, Ottawa, ON, K1S 5B6 Canada}}}

\date{}
\maketitle


\begin{abstract}
\n Some problems of statistics can be reduced to extremal problems of minimizing functionals
of smooth functions defined on the cube $[0,1]^m$, $m\geq 2$.
In this paper, we study a class of extremal problems that is closely connected to
 the problem of testing multivariate independence.  By solving the extremal problem,
we provide a unified approach to establishing weak convergence for a wide class
of empirical processes which emerge in connection with testing independence.
The use of our result is also illustrated by describing the
 domain of local asymptotic optimality of some nonparametric tests of independence.

\bigskip
\n \textbf{Keywords}: boundary-value problem; Green function; multivariate independence;
asymptotic efficiency of test statistics; local asymptotic optimality

\medskip
\n \textbf{AMS Subject Classification}: 62G10; 62G20; 35C05
\end{abstract}

\section{Introduction}

Let $\Xb_i=(X_{i1},\ldots,X_{im}),\;m\geq
2,\,i=1,\ldots,n,$ be independent random vectors with absolutely continuous
distribution function (df) $F$ and marginal df's $F_1,\ldots,F_m$.
An important problem of statistics is to test the hypothesis of (multivariate) independence
\begin{gather}\label{ind}
H_0:F\equiv F_1 F_2\ldots F_m.
\end{gather}
Many test statistics for testing multivariate independence converge weakly toward functionals
of Gaussian random processes under the hypothesis of independence
(see, e.g., \cite{d, nikbook, naznik, shsh}).
If such functionals are regular then the limiting distributions of test statistics are easily derived.
The knowledge of the limiting distribution is important for calculating
(approximate) critical values and finding asymptotic efficiency of the sequence of
 test statistics. Unfortunately,  many independence tests have a complex structure,
 and the corresponding asymptotic theory is not easy to develop. This prevents the use of such tests.
In this paper, by solving a multidimensional extremal problem (Theorem 1),
we prove a result on the weak convergence of empirical processes  (Theorem 2),
 which yields
the limit distribution for a wide class of test statistics for testing independence. Such a class includes,
as extreme cases, the multivariate versions of Cram\'{e}r--von Mises and Blum--Kiefer--Rosenblatt test statistics (see Section 4.1).
 In addition, we illustrate the use of Theorem 1
for describing the domain of local asymptotic optimality of some nonparametric tests of independence (see Section 4.2).

Let $M=\{1,2,\ldots,m\}$ and let $2^M$ be the set of all subsets of $M$.
In this paper, we study a class of extremal problems indexed by a subset $\mathcal{M}$ of $2^M$
that obeys certain restrictions.
For $m=2$ there exist three essentially different choices for $\mathcal{M}$, and the class of extremal problems
is limited. Such problems were studied in \cite[Ch. 5]{nikbook} and \cite{naznik}
in connection with calculating asymptotic efficiency of nonparametric tests of independence in a bivariate setup.
In a general case of $m\geq 2$ the situation is more complicated, and the number of
possible extremal problems for an arbitrary $m$ is not explicitly calculable  (see Remark 1).

The paper
is organized as follows.
In Sections 2 and 3 we formulate and solve the extremal problem of interest.
The solution is obtained by reducing the extremal problem to a non-standard boundary-value problem
and constructing a Green function for the latter.
This part of the paper, with the main result given by Theorem~1, is purely analytical
and might be of independent interest for experts in PDE's.
The rest of the paper is largely statistical, with a focus on various applications of Theorem~1
to the problem of testing independence. In particular, Theorem 1 provides a unified approach to establishing weak convergence for a wide class
of empirical processes with a  multidimensional time parameter which emerge in nonparametric statistics (Theorem 2).

\section{Extremal problem}
Let $I^m=[0,1]^m$ and let $\C_0^m$ be the space of (real-valued) functions that are $m$-times continuously differentiable with respect to
each variable and obey the following boundary conditions:
$${\C}^m_0=\{\Omega\in\C^m(I^m):\Omega(\xb)|_{x_j=0}=0,\;j=1,\ldots,m\},$$
where $\xb=(x_1,\ldots,x_m)\in I^m$.
Define a scalar product on $\C_0^m$ as follows:
\begin{gather}\label{sp}
(\Omega_1,\Omega_2)=\int_{I^m}\omega_1(\xb)\omega_2(\xb)\,{d}\xb,\quad \Omega_1,\Omega_2\in\C_0^m,
\end{gather}
where $\omega_i(\xb)=\dfrac{\partial^m \Omega_i(\xb)}{\partial
x_1\ldots \partial x_m},$ $i=1,2.$ Denote by $\Hb^m$ the closure of
the space $\C^m_0$ under the norm $\|\cdot\|$ induced by the scalar
product (\ref{sp}). For any $m\geq 2$, $\Hb^m$ is a Hilbert space
whose properties are derived similarly to the case $m=2$
studied in \cite{naznik}. More precisely, the following result holds true.

\begin{proposition} (a) The embedding of the space $\Hb^m$ into the space $\C(I^m)$ is compact.
(b) The embedding of the space $\Hb^m$ into the Sobolev space ${\bf W}_2^1(I^m)$ is compact.
\end{proposition}

Part (a) of Proposition 1 implies that any function from
$\Hb^m$ equals zero on any ``left'' side of the cube $I^m$ adjacent
to the origin.

Consider the extremal problem
\begin{gather}
\|\Omega\|_{\Hb^m}\to \min,\label{extpr}\\
 \Omega\in\Hb^m_\mathcal{M},\quad\int_{I^m}\Omega(\xb)\, { d}\mu(\xb)=1,\label{eq12}
 \end{gather}
where $\Hb^m_\mathcal{M}$ is a subset of $\Hb^m$ specified by certain boundary conditions on the ``right'' sides of the cube $I^m$
adjacent to the point ${\bf 1}=(1,1,\ldots,1)$, and $\mu$ is a finite measure on $I^m$.
In order to describe
possible boundary conditions of this extremal problem we need some
notation.

Let $M=\{1,2,\ldots,m\}$ and let $2^M$ be the set of all subsets of $M$.
For any $U\subset M$, denote $\x_U$ the $|U|$-dimensional vector $\xb_U=(x_i:i\in U)$,
where $|U|$ is the cardinality of $U$.
Consider the set $\mathcal{M}\subset 2^M$ such that
\begin{gather}\label{req}
\forall\;U\subset V\subset 2^M,\;\; U\in
\mathcal{M}\;\mbox{ implies}\; V\in\mathcal{M}.
\end{gather}
That is, if set $U$ belongs to
$\mathcal{M}$, then all its ``oversets'' also belong to $\mathcal{M}.$
Define the subset $\Hb^m_\mathcal{M}$ of $\Hb^m$ as follows:
\begin{gather*}
\Hb^m_\mathcal{M}=\{\Omega\in \Hb^m: \Omega(\xb)|_{\xb_U={\bf 1}}=0,\quad U\in\mathcal{M}\}.
\end{gather*}
The
reason for the requirement (\ref{req})  is simple: if $\Omega\in {\Hb^m}$ takes
a zero value on the side $S_U=\{\xb_{U}={\bf 1}\}$, it also takes a zero
value on all the subedges of $S_U$ of less dimension.

For a set $V=(i_1,\ldots,i_l)$ and its complement (in ${M}$)
$V^{c}=(j_1,\ldots,j_k),$ $l+k=m$, put
$$ \partial\xb_V  \partial\xb^2_{V^{c}} =\partial
x_{i_1}\ldots \partial x_{i_l}\partial x^2_{j_1}\ldots \partial
x^2_{j_k}.
$$
By the Lagrange principle rule (see, e.g., \cite{atf}),
the necessary condition of a minimum in (\ref{extpr})--(\ref{eq12}) is reduced to the Euler
equation (in the sense of distributions)
\begin{gather}\label{bp1}
(-1)^m\lambda\frac{\partial^{2m}\Omega(\xb)}{\partial x_1^2\ldots \partial
x_m^2}=\mu(\xb),
\end{gather}
where $\lambda$ is the Lagrange multiplier, and the natural boundary conditions
\begin{gather}\label{bp2}
\frac{\partial^{l+2k}\Omega(\xb)}{\partial{\xb_V} \partial\xb^2_{V^c}}\Big|_{\xb_V =1}=0,\quad\mbox{for any}\;\;
V\notin\mathcal{M},\; V\neq \emptyset.
\end{gather}

\begin{remark} {\rm For any $U\in 2^M$ define an $m$-dimensional vector of Boolean variables $(y_j=\mathbb{I}(j\in U), j=1,\ldots,m)$.
Then $\mathbb{I}(U\in\mathcal{M})$ is a monotone Boolean function
(see \cite{kl}). Denote by $N(m)$ the total number of such functions.
Obviously, the number of the above considered extremal problems is
also equal to $N(m)$. So far, no explicit formula for $N(m)$ as a
function of $m$ has been found. For the asymptotic behavior of $N(m)$,
as $m\to \iy$, see \cite{korsh}.}
\end{remark}

\section{The Green function of the problem}
Recall that Green function of the boundary-value problem (\ref{eq12}), (\ref{bp1}), (\ref{bp2})
is the function ${\cal G}_{\mathcal{M}}(\xb,\xib)$ that satisfies, along with the boundary conditions, the equation (in the sense of distributions)
\begin{gather}\label{Green}
(-1)^m \frac{\partial^{2m} {\cal G}_{\mathcal{M}}(\xb,\xib)}{\partial{x_1^2}\ldots \partial x_m^2}=\delta(\xb-\xib),
\end{gather}
where $\delta(\xb)$ is the Dirac function. It is well known that solution
to the problem (\ref{eq12}), (\ref{bp1}), (\ref{bp2}), and hence solution to the extremal problem
(\ref{extpr})--(\ref{eq12}), can be expressed with the aid of Green function by the formula
\begin{gather}\label{Ome}
\Omega(\xb)=\l^{-1}\int_{I^m}{\cal G}_\mathcal{M}(\xb,\xib)\,{ d}\mu(\xib),\quad \xb\in I^m,
\end{gather}
where the Lagrange multiplier $\lambda$ is found from the integral restriction in (\ref{eq12}) and
has the form
\begin{gather}\label{lam}
\l=\iint_{I^m\times I^m}{\cal G}_{\mathcal{M}}(\xb,\xib)\,{ d}\mu(\xb)\,{d}\mu(\xib).
\end{gather}

For $\xb, \xib\in I^m$ and for $V\subset M$ as before, define the functions
$$k_{V}(\xb,\xib)=k_{i_1}(\xb,\xib)\ldots k_{i_l}(\xb,\xib),\quad
K_{V^c}(\xb,\xib)=K_{j_1}(\xb,\xib)\ldots K_{j_k}(\xb,\xib), $$
where
\begin{gather}\label{Kk}
k_j(\xb,\xib)=x_j \xi_j,\quad K_j(\x,\xib)=\min(x_j, \xi_j ),\quad j=1,\ldots,m.
\end{gather}

The first main result of the paper is stated as follows.
\begin{theorem} The Green function of the boundary-value problem {\rm(\ref{eq12})}, {\rm (\ref{bp1})}, {\rm(\ref{bp2})} is
\begin{gather}\label{Gf}
{\cal G}_{\mathcal{M}}(\xb,\xib)=K_M(\xb,\xib)-\sum_{U\in\mathcal{M}}a_U K_{U^c}(\xb,\xib) k_{U}(\xb,\xib),
\end{gather}
with the coefficients $a_U$ defined recurrently by
\begin{gather}\label{rec}
\sum_{V\subset U\atop V\in{\mathcal{M}}}a_V=1,\quad\mbox{for all}\;\; U\in\mathcal{M}.
\end{gather}
\end{theorem}

\begin{proof} First, note that
\begin{gather*}
-\frac{\partial^2 K_j(\xb,\xib)}{\partial{x_j^2}}=\delta(x_j-\xi_j),\quad \frac{\partial^2 k_j(\xb,\xib)}{\partial{x_j^2}}=0.
\end{gather*}
Therefore the function ${\cal G}_{\mathcal{M}}$
satisfies (\ref{Green}) with an arbitrary choice of the constants $a_U$.

Further, by (\ref{req}) for any nonempty set $V=(i_1,\ldots,i_l)\notin{\mathcal{ M}}$ and
for any $U\in\mathcal M$ there exists $j\in U^c\cap V$, and hence
$$\frac {\partial^l K_{U^c} ({\mathbf x}, \mbox{\boldmath$\xi$})
k_U ({\mathbf x}, \mbox{\boldmath$\xi$})} {\partial {\mathbf
x}_V}\Big|_{{\mathbf x}_V={\bf 1}}= \frac {\partial K_j({\mathbf x},
\mbox{\boldmath$\xi$})} {\partial
x_j}\Big|_{x_j=1}\times\dots=0.$$
Thus all summands in (\ref{Gf}) satisfy (\ref{bp2}).

In order to prove that the function ${\cal G}_{\mathcal M}$ vanishes on $\{\mathbf{
x}_U={\bf 1}\}$  for some $U\in\mathcal M$, we represent $\mathcal M$ as the
union $\mathcal M=\bigcup\limits_{V\supset U}{\mathfrak M}_V$ of the disjoint sets
$${\mathfrak M}_V=\{W\in\mathcal M:\ W\cup U=V\}$$
(in the definition of ${\mathfrak M}_V$ the basic property (\ref{req}) is used).

The key observation is that for any $W\in{\mathfrak M}_V$,
$$K_{W^c} ({\mathbf x}, \mbox{\boldmath$\xi$})k_W ({\mathbf x},
\mbox{\boldmath$\xi$}) \big|_{{\mathbf x}_U={\bf 1}}= K_{V^c} ({\mathbf
x}, \mbox{\boldmath$\xi$})k_V ({\mathbf x}, \mbox{\boldmath$\xi$})
\big|_{{\mathbf x}_U={\bf 1}}.$$
Therefore formula (\ref{Gf}) gives
\begin{multline*}
{\cal G}_{\mathcal M}({\mathbf x},
\mbox{\boldmath$\xi$})\big|_{{\mathbf x}_U={\bf 1}}= K_{U^c} ({\mathbf
x}, \mbox{\boldmath$\xi$})k_U ({\mathbf x}, \mbox{\boldmath$\xi$})
\big|_{{\mathbf x}_U={\bf 1}}
\cdot\Big(1-\sum\limits_{W\in{\mathfrak M}_U}a_W\Big)-\\
-\sum\limits_{V\supset U\atop V\neq U} K_{V^c} ({\mathbf x},
\mbox{\boldmath$\xi$})k_V ({\mathbf x}, \mbox{\boldmath$\xi$})
\big|_{{\mathbf x}_U={\bf 1}} \cdot\sum\limits_{W\in{\mathfrak M}_V}a_W.
\end{multline*}
Due to (\ref{rec}) the coefficients on the right-hand side vanish, and
${\cal G}_{\mathcal M}({\mathbf x},
\mbox{\boldmath$\xi$})\big|_{{\mathbf x}_U={\bf 1}}=0$. This
completes the proof.
\end{proof}

\begin{remark} {\rm  Return to the extremal problem (\ref{extpr})--(\ref{eq12}) and consider the following three sets of boundary conditions: (a) there
are no restrictions on $\Omega\in{\Hb}^m$ except for those that
specify the space ${\Hb}^m$, (b) $\Omega\in{\Hb}^m$ equals zero on
any $(m-1)$-dimensional side of $I^m$, and (c) $\Omega\in{\Hb}^m$
equals zero at the point ${\bf 1}=(1,\ldots, 1)$. Then $\mathcal{M}=\emptyset$, $\mathcal{M}=2^M$, $\mathcal{M}=\{M\}$, respectively,
and by Theorem~1 the corresponding Green functions are  $\prod_{j=1}^m K_j(\xb,\xib)$, $\prod_{j=1}^m (K_j(\xb,\xib)-k_j(\xb,\xib))$,
 $\prod_{j=1}^m K_j(\xb,\xib)-\prod_{j=1}^m k_j(\xb,\xib)$.
These are covariance functions of the classical Gaussian random
processes. They correspond to an $m$-dimensional Brownian sheet, an $m$-dimensional Brownian pillow, and
an $m$-dimensional tucked Brownian sheet, respectively.
The latter two arise as limiting
processes in nonparametric testing of independence (see Section 4.1 for details).}
\end{remark}

\section{Connection to testing independence}
This section illustrates the use of Theorem 1 for some efficiency issues that
emerge in the problem of  testing multivariate independence.
We use a general dependence model which is rather popular in the present context.

Let $\Xb_i=(X_{i1},\ldots,X_{im}),\;m\geq
2,\,i=1,\ldots,n,$ be independent random vectors with absolutely continuous
df $F$ and marginal df's $F_1,\ldots,F_m$.  Consider the problem of testing the hypothesis of independence
(\ref{ind}).

For the absolutely continuous df $F$, a \textit{copula} $\phi$ of $F$ is a (unique) df with uniform univariate margins, such that
$\phi(F_1(x_1),\ldots,F_m(x_m))=F(x_1,\ldots,x_m)$.
When testing for independence using distribution-free (independent of $F_1,\ldots,F_m$)
test statistics, one can assume uniform margins and define the model in terms of copulas.

Let $\{F_\t:\t\geq 0\}$ be the family of
absolutely continuous copulas  $F_\t$ that are monotone in $\t$ (the case
$\t=0$ corresponds to independence) and that satisfy the following common regularity conditions, cf. \cite{genest, quessy}:
\begin{enumerate}
\item[{\bf(A1)}] the density $f_\t(\xb)=\partial^m F_\t(\xb)/\partial x_1\ldots\partial x_m$
admits a square-integrable derivative $\dot{f}_0(\xb)$ of
$f_\t(\xb)$ with respect to $\t$ at $\t=0$ for every $\xb=(x_1,\ldots,x_m)\in I^m$,
and the function $\sqrt{f_{\t}(\xb)}$ is differentiable in quadratic
mean at $\t=0$, i.e.,
$$\int_{I^m}\left(\sqrt{f_\t(\xb)}-1-\frac12 \t \dot{f}_0(\xb)\right)^2\,{d}\xb=o(\t^2),\quad \t\to 0,$$

\item[{\bf(A2)}] $F_\t(\xb)$ is differentiable with respect to $\t$ in a small neighborhood of $\t=0$ for every $\xb\in I^m,$ and
the following identity holds for every $\xb\in I^m:$
$$\dot{F}_0(\xb)=\lim_{\t\to 0}\frac{\partial}{\partial\t}F_{\t}(\xb) =\int_0^{x_1}\ldots\int_0^{x_m}\dot{f}_0(y_1,\ldots,y_m)\,{\rm d} y_1\ldots\,{d} y_m;$$
and, in addition,
$$\int_{I^m}\dot{F}_{0}(\xb) \,{d}\xb=\lim_{\t\to
0}\frac{\partial}{\partial\t}\int_{I^m}F_{\t}(\xb)\,{d}\xb.$$
\end{enumerate}

The function $\dot{F}_0(\xb)$ is sometimes called the \textit{dependence function}.
Assumption {\bf (A2)} implies that $F_\t(\xb)$ can be written in
terms of $\dot{F}_0$ as follows:
\begin{gather*}
F_\t(\xb)=\prod_{j=1}^m x_j +\t \dot{F}_0(\xb)+o(\t),\quad\xb\in I^m,\quad \t\to 0,
\end{gather*}
where, by the property of a
multivariate copula,
the boundary conditions
\begin{gather}\label{nerav3}
\dot{F}_0(\xb)|_{x_k=0}=0,\quad \quad \dot{F}_0(\xb)|_{\xb_U={\bf 1}}=0,\quad \mbox{for any}\;\; U\in 2^M,\;\;|U|=m-1
\end{gather}
are satisfied. As an alternative to $H_0:\t=0$ we consider the hypothesis $H_1:\t>0$.
In what follows, the underlying df is assumed to belong to the family $\{F_\t:\t\geq 0\}$.

The family $\{F_\t:\t\geq 0\}$ produces a sequence of locally asymptotically normal experiments.
Indeed, let $\Pb_\theta$ be the probability distribution used to calculate the df $F_\theta$.
Then the full observation is a single observation
from the product $\Pb^{n}_{\theta}$ of $n$ copies of $\Pb_\theta$.
In view of condition {\bf (A1)}, the sequence of statistical experiments $\{\Pb^n_{h/\sqrt{n}}:h\geq 0\}$, indexed by a local parameter $h=\sqrt{n}\t$, is
locally asymptotically normal (LAN) at the point $h=0$, that is,
\begin{equation*}
\log\frac{d\Pb^n_{h/\sqrt{n}}}{d\Pb^n_0}(\Xb_1,\ldots,\Xb_n)=
h\Delta_{n,0}-\frac12 h^{2}I_0+o_{\Pb^n_0}(1), \quad n\to \iy,
\end{equation*}
where $I_{\theta}=\int_{I^m}{\dot{f}^2_\theta({\bf
x})}/{f_\theta(\xb)}\,{\rm d}\xb$ is the Fisher
information in the family $\{F_\theta,\, \theta \geq
0\}$
and $\Delta_{n,0}=n^{-1/2}\sum_{i=1}^n
\dot{f}_0(\Xb_i)\stackrel{\P^{n}_{0}}{\rightsquigarrow}{
N}(0,I_{0})$.
The symbol $\stackrel{\P^{n}_{0}}{\rightsquigarrow}$ denotes convergence in $\P^{n}_{0}$-distribution.
The LAN property ensures the mutual contiguity of the sequences of distributions $\{\Pb^n_{h/\sqrt{n}}\}$ and $\{\Pb^n_{0}\}$, and
facilitates the calculation of Pitman efficiency of asymptotically normal test statistics (see \cite[Ths. 14.7, 15.4]{vaart}).

\subsection{Asymptotic efficiency of independence tests}
Consider testing the hypothesis of independence using distribution-free statistics of $\Xb_1,\ldots, \Xb_n$
in the case when
some of the $F_j$'s, $j=1,\ldots,m,$ are \textit{known}, while the others are \textit{unknown}.
Denote by
$\mathbb{F}_n$ the multivariate empirical df that corresponds to $F$, and
denote by $\mathbb{F}_{j,n}$ the marginal empirical df based on
$X_{1j},\ldots,X_{nj}$, $j=1,\ldots,m.$
For a set $V=(i_1,\ldots,i_l)$ and its complement (in
${M}$) $V^{c}=(j_1,\ldots,j_k),$ $l+k=m$,
put
\begin{gather*}
F_{V}(\xb)=F_{i_1}(x_{i_1})\ldots F_{i_l}(x_{i_l}),\quad
{ d} F_{V}(\xb)={ d} F_{i_1}(x_{i_1})\ldots { d}F_{i_l}(x_{i_l}),\\
\mathbb{F}_{V^c,n}(\xb)=\mathbb{F}_{j_1,n}(x_{j_1})\ldots \mathbb{F}_{j_k,n}(x_{j_k}),
\quad { d} \mathbb{F}_{V^c,n}(\xb)={d} \mathbb{F}_{j_1,n}(x_{j_1})\ldots {d} \mathbb{F}_{j_k,n}(x_{j_k}).
\end{gather*}
For $p=1,2,\ldots,$ consider the class $\{B_{V,n}^p:V\subset M\}$ of test statistics, cf. \cite{dug},
\begin{gather*}
B^p_{V,n}=\int_{\Rbb^m}(\mathbb{F}_n(\xb)-F_{V}(\xb) \mathbb{F}_{V^c,n}(\xb))^p  \,{d}F_{V}(\xb)\,{ d}\mathbb{F}_{V^c,n}(\xb).
\end{gather*}
Suppose that the margins $F_{j}$, $j\in V$, are known. Then
the test statistics $B^p_{V,n}$, $V\subset M$, $p=1,2,\ldots,$ are distribution-free under the null hypothesis, and
the study of their behaviour under $H_0$
 can be done when $F(\xb)$ is a uniform distribution on $I^m$, which will be assumed from now on.
Choosing $V=M$
yields the Cram\'{e}r--von Mises-type statistics
\begin{eqnarray*}
B^p_{M,n}&=&\int_{I^m}(\mathbb{F}_n(\xb)-\xb_M)^p \, { d} \xb_M,
\end{eqnarray*}
where $\xb_M=x_{1}\ldots x_{m}$ and ${\rm d}\xb_M= {\rm d}\xb={\rm d}x_{1}\ldots {\rm d}x_{m}.$
On the other hand, setting $V=\emptyset$ leads to the Blum--Kiefer--Rosenblatt-type statistics
\begin{eqnarray*}
B_{\emptyset,n}^p&=&\int_{I^m}(\mathbb{F}_n(\xb)- \mathbb{F}_{M,n}(\xb))^p
 \, {d}\mathbb{F}_{M,n}(\xb).
 \end{eqnarray*}

An important step in calculating asymptotic efficiency of the sequence of  test statistics $\{B^p_{V,n}\}_{n\geq 1}$ lies in
showing the weak convergence of
$$W_{V,n}(\xb)=\sqrt{n}\left(\mathbb{F}_n(\xb)-\xb_{V} \mathbb{F}_{V^c,n}(\xb)\right),\quad \xb\in I^m,$$
to a limiting Gaussian process under the null hypothesis. This is achieved with the aid of Theorem 1, and is stated below as Theorem 2.
The Skorohod space $D(I^m)$ that appears in the statement of Theorem 2 and  generalizes the well-known space $D[0,1]$,
is described in \cite{neuh}, where some of its properties are also derived.

\begin{theorem} Assume that $F(\xb)$ is a uniform distribution on $I^m$.
Then for any $V~\subset~M$ the empirical process $W_{V,n}(\xb)$ converges weakly in the Skorohod space $D(I^m)$ to
a centered Gaussian process $W_V(\xb)$ with covariance function
${\cal G}_{\mathcal{M}_V}(\xb,\xib)$ given by {\rm (\ref{Gf})}, where the set ${\mathcal M}_V\subset 2^M$ is defined as follows:
\begin{gather}\label{MMM}
{\mathcal{M}_V}=\left\{M, M\setminus U: U\subset V^c, |U|= 1 \right\}.
\end{gather}
\end{theorem}

\begin{proof}
The proof follows the pattern of \cite[Sec. 4, 5]{neuh}, where, among others, weak convergence of the process
$W_{M,n}(\xb)=\sqrt{n}\left(\mathbb{F}_n(\xb)-\xb_{M}\right)$
 to a tucked Brownian sheet is established.
Therefore, most details on the convergence of finite-dimensional distributions of $W_{V,n}(\xb)$ to a multivariate normal distribution,
and the proof of tightness of the family of associated probability measures are omitted.
A pivotal point of our proof is obtaining
an expression for the covariance of the limiting process.
This is accomplished by appealing to Theorem 1.

For some $j\notin V$ let $U=M\setminus\{j\}$. Then $\mathbb{F}_n(\xb)|_{x_U={\bf 1}}=\mathbb{F}_{j,n}(x_j)$, so that the process $W_{V,n}(\xb)$
is pinned down to zero on $\{\xb_U={\bf1}\}$ for any $U\in \mathcal{M}$ such that $V\subset U$ and $|U|=m-1. $

Now consider the boundary-value problem {\rm(\ref{eq12})}, {\rm (\ref{bp1})}, {\rm(\ref{bp2})} with the set $\mathcal{M}={\mathcal M}_V$
given by (\ref{MMM}). The respective boundary condition takes the form
$$\Omega\in{\Hb}^m,\quad \Omega(\xb)|_{\xb_U={\bf1}}=0\;\;\mbox{for any}\;  U\in \mathcal{M}\;\;\mbox{such that}\; \;V\subset U,\;|U|=m-1.   $$
The function $\Omega(\xb)$ equals zero exactly on those sides $\{\xb_U={\bf 1}\}$ of the cube $I^m$ where
the empirical process $W_{V,n}(\xb)$ vanishes.
This observation together with Theorem 1 gives the required expression for the covariance function.
The proof is completed.
\end{proof}

In the two extreme cases, when (1) $V=M$ and (2) $V=\emptyset$, Theorem 2 yields a well-known result (see, e.g., \cite{blum, neuh, d}).
Indeed, in the first case the set in (\ref{MMM}) reduces to $\mathcal{M}_M=\{M\}$ and the covariance function of $W_M(\xb)$  is
\begin{gather}
{\cal G}_{\mathcal{M}_M}(\xb,\xib)=\prod_{j=1}^m K_j(\xb,\xib)-\prod_{j=1}^m k_j(\xb,\xib),\label{m}
\end{gather}
The respective boundary condition
has the form
$\Omega\in{\Hb}^m,\;\Omega(\xb)|_{\xb_M={\bf 1}}=0$.
In the second case (\ref{MMM}) becomes
$$\mathcal{M}_\emptyset=\{M,M\setminus\{1\},M\setminus\{2\},\ldots,M \setminus\{m\}\},$$ and
the covariance function of $W_\emptyset(\xb)$
is
\begin{gather}
{\cal G}_{\mathcal{M}_\emptyset}(\xb,\xib)=
\prod_{j=1}^m K_j(\xb,\xib)-
\sum_{j=1}^m K_j(\xb,\xib)\prod_{i\neq j}k_i(\xb,\xib)+(m-1)\prod_{j=1}^m k_j(\xb,\xib)\label{mm}.
\end{gather}
The respective boundary condition
takes the form
$\Omega\in{\Hb}^m,\; \Omega(\xb)|_{\xb_U=1}=0$ for any $U\in 2^M$ such that $|U|=m-1$.

\medskip
For large sample sizes, the quality of test statistics $B_{V,n}^p$, $V\subset M$, $p\geq 1$, can be judged by looking at their Bahadur efficiency.
This kind of asymptotic efficiency is quantified by the \textit{Bahadur exact slope}.
Finding the Bahadur exact slope of a sequence of test statistics requires
the law of large numbers under the alternative, and the rough large deviation asymptotics under the null hypothesis.
The problem of large deviation asymptotics consists in minimizing the Kullback--Leibler information
over a subset of distribution functions that depends on the structure of the test statistic.
When $p=1$ the above minimization problem is reduced, by using variational methods, to the boundary-value problem (\ref{eq12}), (\ref{bp1}), (\ref{bp2})
with $\mu$ being the Lebesgue measure on $I^m$ and $\mathcal{M}=\mathcal{M}_V$ given by (\ref{MMM}),
whose solution provides the main contribution to the initial problem (see \cite[Ch. 5]{nikbook} for details).

For example, finding the rough large deviation
asymptotics of the statistic
$$B_{\emptyset,n}^1=\int_{I^m}(\mathbb{F}_n(\xb)-\mathbb{F}_{M,n}(\xb)) \,{ d}\mathbb{F}_{M,n}(\xb)$$
is largely reduced to the boundary-value problem {\rm(\ref{eq12})}, {\rm (\ref{bp1})}, {\rm(\ref{bp2})}
with $\mathcal{M}= \mathcal{M}_\emptyset$
and for sufficiently small $t>0$, cf. formula (5.3.29) of \cite{nikbook},
\begin{gather}\label{bo}
\lim_{n\to \iy}n^{-1}\log \P_{H_0}(B_{\emptyset,n}^1\geq t)=-\frac12 \l_0 t^2+\sum_{j\geq 3} c_j t^j,
\end{gather}
where the series on the right-hand side is convergent and, cf. (\ref{Ome}) and (\ref{lam}),
\begin{gather*}
\l_0=\left(\iint_{I^m\times I^m}{\cal G}_{\mathcal{M}_\emptyset}(\xb,\xib)\,{d}\xb\,{d}\xib\right)^{-1}
=\frac{4^m}{(4/3)^m-m/3 -1}
\end{gather*}
(see \cite[Sec. 5.3]{nikbook}).
Since, under the alternative,
\begin{gather*}
B_{\emptyset,n}^1 \stackrel{\P_\t}{\rightarrow} \t\int_{I^m}\dot{F}_0(\xb)\,{ d}\xb,\quad n\to \iy,
\end{gather*}
where the symbol $\stackrel{\P^{n}_{0}}{\rightarrow}$ denotes convergence in $\P^{n}_{0}$-probability,
it follows from (\ref{bo}) and Theorem 1.2.2 of \cite{nikbook} that the Bahadur exact slope $c_{B^1_\emptyset}(\t)$
of the sequence $\{B_{\emptyset,n}^1\}_{n\geq 1}$ satisfies as $\t \to 0$
\begin{gather}\label{slope}
c_{B^1_{\emptyset}}(\t)\sim \t^2 \frac{4^m}{(4/3)^m-m/3 -1}\left(\int_{I^m}\dot{F}_0(\xb) \,{d}\xb\right)^2.
\end{gather}

Similarly, in a general case, for the model under consideration a routine computation leads to the following result.
\begin{proposition} For an arbitrary  $V\subset M$ the Bahadur exact slope of the sequence $\{B_{V,n}^1\}_{n\geq 1}$ satisfies
$$c_{B^1_{V}}(\t)\sim \t^2\left(\iint_{I^m\times I^m}{\cal G}_{\mathcal{M}_V}(\xb,\xib)\,{d}\xb\,{ d}\xib \right)^{-1}\left(\int_{I^m}\dot{F}_0(\xb)\,{ d}\xb\right)^2,\quad \t\to 0,$$
where the function ${\cal G}_{\mathcal{M}_V}(\xb,\xib)$ is the same as in Theorem 2.
\end{proposition}

Compared to $B_{V,n}^1$, the evaluation of the local Bahadur efficiency of the tests based on $B_{V,n}^p$, $p\geq 2$,
is much more complicated. For example, when $p=2$, the efficiency problem is reduced to
calculating the principal eigenvalue of the integral operator with kernel
${\cal G}_{\mathcal{M}_V}$ (see \cite[Ch. 5]{nikbook} for details).

\medskip

A more complex empirical process that vanishes completely at the boundary of $I^m$ was studied, for example, in \cite{neuh} and \cite{d1}.
Such a process converges weakly to the $m$-dimensional Brownian pillow with covariance function
 \begin{gather}\label{cf}
{{\cal G}}_{ 2^M }(\xb,\xib)=\prod_{j=1}^m \left(K_j(\xb,\xib)-k_j(\xb,\xib)\right),\quad \xb,\xib\in I^m,
\end{gather}
which corresponds to the Green function (\ref{Gf}) with $\mathcal{M}=2^M$ (see Remark 2),
and appears in connection with testing multivariate independence in the following context.

Consider $W_{\emptyset,n}(\xb)=\sqrt{n}(\mathbb{F}_n(\xb)-\prod_{j=1}^m \mathbb{F}_{j,n}(x_j))$.
The corresponding limiting process $W_\emptyset(\xb)$ has the covariance function
${\cal G}_{\mathcal{M}_\emptyset}(\xb,\xib)$ which  coincides with ${\cal G}_{2^{M}}(\xb,\xib)$ when $m=2$.
However, for $m\geq 3$ the situation changes. The process $W_\emptyset(\xb)$ does
not vanish completely at the facets of $I^m$ adjacent to the point $\xb={\bf 1}$,
but  only does so at the one-dimensional edges (see (\ref{mm})).
This ``disappointing'' property is overcome by the \textit{tied-down} empirical process (see, e.g., \cite{neuh, d1})
\begin{gather*}\hat{W}_{\emptyset,n}(\xb)=\sqrt{n}\left(\mathbb{F}_n(\xb)-\sum_{k=1}^m (-1)^{k-1}
\sum_{U\subset M:|U|=k}\xb_U\cdot \mathbb{F}_n(\xb)|_{\xb_U={\bf 1}}\right),\quad \xb\in I^m,
\end{gather*}
which can be equivalently written in the form
\begin{gather*}
\hat{W}_{\emptyset,n}(\xb)=\frac{1}{\sqrt{n}}\sum_{i=1}^n \prod_{j=1}^m \left(\mathbb{I}(X_{ij}\leq x_j)-x_j\right),\quad \xb\in I^m.
\end{gather*}
Under the null hypothesis, $\hat{W}_{\emptyset,n}(\xb)$ converges weakly in the Skorokhod space $D(I^m)$ toward
an $m$-dimensional Brownian pillow $\hat{W}_\emptyset(\x)$ with covariance function (\ref{cf}).
This fact seems to be established for the first time in \cite{neuh}.
The corresponding test statistics (for $m\geq 3$) take the form,  cf. \cite{dug,d},
\begin{eqnarray*}
\hat{B}^p_{m,n}&=&\int_{I^m}\left(\mathbb{F}_n(\xb)-
\sum_{k=1}^m (-1)^{k-1} \sum_{U\subset M:|U|=k}\xb_U\cdot \mathbb{F}_n(\xb)|_{\xb_U={\bf 1}} \right)^p\,{ d}\xb.
\end{eqnarray*}
For $p=1$,
under the hypothesis of independence, the limiting distribution of $\sqrt{n} \hat{B}^1_{m,n}$ is normal
with zero mean and variance, cf. Th. 6 of \cite{d},
\begin{eqnarray*}
\hat{\s}_m^2(0)&=&
\iint_{I^m\times I^m} {{\cal G}}_{2^{M}}(\xb,\xib)\,{ d}\xb\,{ d}\xib=12^{-m}.
\end{eqnarray*}
Therefore, by  Le Cam's third lemma,  for all $h\geq 0$,
$$\frac{\sqrt{n}(\hat{B}^1_{m,n}-\hat{\mu}_m(h/\sqrt{n}))}{\hat{\s}_m(0)} \stackrel{\P^{n}_{h/\sqrt{n}}}{\rightsquigarrow}{\cal N}(0,1),\quad n\to \iy,$$
where, using (\ref{nerav3}),
\begin{gather*}\hat{\mu}_m(\theta) =\t\int_{I^m}\left(\dot{F}_0(\xb)-
\sum_{k=1}^{m-2} (-1)^{k-1}\sum_{U\subset M: |U|=k}\xb_U\cdot \dot{F}_0(\xb)|_{\xb_U={\bf 1}}\right) \,d\xb
\end{gather*}
 The symbol $\stackrel{\P^{n}_{\t}}{\rightsquigarrow}$ denotes convergence in $\P^{n}_{\t}$-distribution of a random sample
drawn from $F_\t$. In view of Theorem 14.7 in \cite{vaart}, the squared Pitman slope of the sequence $\{\hat{B}_{m,n}^1\}_{n\geq 1}$
 is $$\left(\frac{\hat{\mu}^{\prime}_m(0)}{\hat{\s}_m(0)}\right)^2=12^{m}\left(\int_{I^m}\left(\dot{F}_0(\xb)-
\sum_{k=1}^{m-2} (-1)^{k-1}\sum_{U\subset M: |U|=k}\xb_U\cdot \dot{F}_0(\xb)|_{\xb_U={\bf 1}}\right)\,d\xb\right)^2.$$

\begin{remark} {\rm Similar to $\hat{W}_{\emptyset,n}(\xb)$ one can get other empirical processes with a limiting covariance
of the form (\ref{Gf}) by subtracting from $\mathbb{F}_n$ linear combinations of empirical processes of less dimension.}
\end{remark}

\subsection{Local asymptotic optimality of independence tests}
An interesting statistical problem that leads to the extremal problem (\ref{extpr})--(\ref{eq12})
is that of local asymptotic optimality of independence tests.

Consider testing the hypothesis of independence
$H_0:\;\theta=0$ versus the alternative $H_1:\theta>0.$
Two commonly used measures for judging the quality of testing are
the Pitman slope and the Bahadur local index.
Under both approaches, the measure of efficiency 
of a given test statistic $T_n=T(\Xb_1,\ldots,\Xb_n)$
has an upper bound (see, e.g., \cite[Th. 15.4]{vaart} and \cite{Bah}), which yields the inequality
\begin{gather}\label{ub}
b_T(\dot{F}_0)\leq
\int_{I^m}\dot{f}_0^2(\xb)\,{d}\xb.
\end{gather}
Here $b_T$ is a homogeneous functional of degree 2 defined on the space ${\bf H}^m$ that depends on a structure of  $T_n$
and measures efficiency of the corresponding test.
For the Bahadur efficiency the upper bound (\ref{ub})  is a local version of the Bahadur--Ragavachari inequality (see \cite{Bah}).
For the test based on $T_n$, the closer $b_T(\dot{F}_0)$ is to $\int_{I^m}\dot{f}_0^2(\xb)\,d\xb$,
the better the family $\{F_\t:\t\geq 0\}$ is.
Thus, in order to describe the domain of Bahadur and/or Pitman optimality of the sequence of test statistics $\{T_n\}$,
we need to know for which dependence function $ \dot{F}_0$ equality in  (\ref{ub}) is attained.
If $b_T$ is the square of a linear functional, this leads to extremal problem (\ref{extpr})--(\ref{eq12}).
For $m=2$ some applications related to establishing Bahadur optimality of independence tests
can be found in \cite{naznik}. Here we cite two examples from \cite{naznik} with \textit{non-Lebesgue} measure $\mu(\xb)$
in the problem (\ref{extpr})--(\ref{eq12}) that
corresponds to the integration over diagonal(s) of $I^m$. These examples are connected to testing independence using the Gini
rank statistic and Spearman's footrule.

Let $\Xb_i=(X_{i1},\ldots,X_{im}),\;m\geq
2,\,i=1,\ldots,n,$ be as before.
 Denote by $R_{ij}$ the rank of $X_{ij}$ among
$X_{1j},\ldots, X_{nj},$ $i=1,\ldots,n,$ $j=1,\ldots, m.$
Recall that the Gini rank coefficient is defined for $m=2$ by $$r_G=\frac{2}{D_n}\sum_{i=1}^n \left(|n+1-R_{i1}-R_{i2}|-|R_{i1}-R_{i2}|\right),$$
where $D_n=n^2$ if $n$ is even and $D_n=n^2-1$ if $n$ is odd, and inequality (\ref{ub}) takes the form (see \cite{naznik})
\begin{gather}\label{gini}
24\left(\int_0^1 \left(\dot{F}_0 (x,x)+\dot{F}_0 (1-x,x)\right)\, { d}x\right)^2\leq \int_{I^2}\dot{f}_0^2 (\xb)\, { d}\xb.
\end{gather}
In this case $\mu(\xb)=\delta(x_1-x_2)+\delta(1-x_1-x_2)$, $\xb=(x_1,x_2)\in I^2$,
and equality in (\ref{ub})
is attained for the function, cf. Theorem 1,
\begin{gather}\label{opt1}
\dot{F}_0(\xb)=C \int_{I^2} {\cal G}_{\mathcal{M}_\emptyset}(\xb,\xib)\,{ d}\mu(\xib),
\quad C>0,
\end{gather}
where, in view of (\ref{nerav3}), ${\cal G}_{\mathcal{M}_\emptyset}(\xb,\xib)$ is given by (\ref{mm}).
Integrating in (\ref{opt1}) yields
$$\dot{F}_0(\xb)=C\left(|x_1-x_2|^3-|x_1+x_2-1|^3-3(x_1^2+x_2^2)+3(x_1+x_2-1)\right),\quad C>0.$$

For the Spearman footrule based on the statistic
$r_f=\sum_{i=1}^n |R_{i1}-R_{i2}|$
the local Bahadur index on the left-hand side of (\ref{ub}) equals
$$b_{r_f}(\dot{F}_0)=90\left(\int_0^1\dot{F}_0(x,x)\,{ d}x\right)^2,$$ which corresponds to the measure
$\mu(\xb)=\delta(x_1-x_2)$, $\xb=(x_1,x_2)\in I^2$. Therefore the optimal dependence function has the form
\begin{eqnarray*}
\dot{F}_0(\xb)&=&C \int_{I^2} {\cal G}_{\mathcal{M}_\emptyset}(\xb,\xib)\,{ d}\mu(\xib)\\
&=& C\left(|x_1-x_2|^3-(x_1+x_2)^3+2x_1 x_2(x_1^2 +x_2^2 +2)\right),\quad C>0.
\end{eqnarray*}

Another interesting application, when $m\geq 2$, is connected to Pitman optimality of a multivariate Spearman's rho, cf. \cite{shsh, quessy}:
\begin{eqnarray*}
S_{m,n}&=&\frac{1}{C_m}\left\{{n}^{-1}\sum_{i=1}^n \prod_{k=1}^m (n+1-R_{ij}) -
{\left(\frac{n+1}{2}\right)}^m\right\},
\end{eqnarray*}
where $C_m={n}^{-1}\sum_{i=1}^n i^m - {\left({(n+1)}/{2}\right)}^m$
is a normalizing factor.
The statistic $S_{m,n}$ is a sample counterpart of the functional
$s_m(F)=\frac{2^m(m+1)}{2^m-m-1}\left(\int F \,{ d} F_1\ldots \,{d}F_m
-{2^{-m}}\right)$.

As before, consider the family $\{F_\t:\t\geq 0\}$ of
absolutely continuous copulas
that satisfy {\bf (A1)} and {\bf (A2)}.
Under the Pitman approach, $\t=\t_n=h/\sqrt{n}$, where $h\geq 0$ is a {local parameter},
and for all $h\geq 0$, cf. \cite[Cor. 1]{quessy},
$$\frac{\sqrt{n}(S_{m,n}-\mu_m(h/\sqrt{n}))}{\s_m(0)} \stackrel{\P^{n}_{h/\sqrt{n}}}{\rightsquigarrow}{\cal N}(0,1),\quad n\to \iy,$$
where
$\mu_m(\theta)=\frac{2^m(m+1)}{2^m-m-1}\,\theta\int_{I^m}\dot{F}_0(\xb)\,
{\rm d}\xb$ and $ \s^2_m(0)=\frac{(m+1)^2\left((4/3)^m-m/3-1\right)}{(2^m
 -m-1)^2}.$
In view of Theorem 15.4 in \cite{vaart}, inequality (\ref{ub}) takes the form
\begin{gather}\label{bound}
\frac{4^m}{(4/3)^m-m/3-1}\left(\int_{I^m}\dot{F}_0(\xb)\,{ d}\xb\right)^2\leq \int_{I^m}\dot{f}_0^2(\xb)\,{d}\xb.
\end{gather}
Then, the application of Theorem 1 yields that the sequence of test statistics $\{S_{m,n}\}_{n\geq 1}$
is Pitman optimal if and only if, cf . \cite[Ths. 2, 3]{step} and \cite[Sec. 4.4]{quessy},
\begin{gather}\label{op}
\dot{F}_0(\xb)=C\prod_{j=1}^m x_j \left(\prod_{j=1}^m (2-x_j)+\sum_{j=1}^m x_j -(m+1)\right),\quad
\xb\in I^m,\quad C>0.
\end{gather}
Indeed, the test based on $S_{m,n}$ is the ``best''
for those dependence functions $\dot{F}_0$ that deliver equality in inequality (\ref{bound}).
Thus, taking into account (\ref{nerav3}), we minimize the functional $\int_{I^m}\dot{f}_0^2(\xb)\,{d}\xb$
on the space $\Hb^m$ subject to
\begin{gather*}
\int_{I^m}\dot{F}_0(\xb)\,{d}\xb =1,\quad \dot{F}_0(\xb)|_{\xb_U={\bf 1}}=0,\quad \mbox{for any}\;\; U\in 2^M,\;\;|U|=m-1.
\end{gather*}
By the results of Sections 2 and 3, including Theorem 1, the functional $\int_{I^m}\dot{f}_0^2(\xb)\, { d}\xb$ is minimized when
\begin{gather}\label{opt}
\dot{F}_0(\xb)=\l^{-1}\int_{I^m}{\cal G}_{\mathcal{M}_\emptyset}(\xb,\xib)\, { d}\xib,
\end{gather}
where ${\cal G}_{\mathcal{M}_\emptyset}$ is given by (\ref{mm}).
 By homogeneity of inequality (\ref{bound})
the extremal function is defined up to a positive constant.
Integrating in (\ref{opt}) yields (\ref{op}). 

\begin{remark} {\rm For $m=2$ the test statistics $B_{\emptyset,n}^1$ and $S_{m,n}$ are known to be asymptotically
equivalent (see, e.g., \cite[Ch. 5]{nikbook}). The results of this section  extend this property to all $m\geq 2$.
Indeed, due to (\ref{slope}) and (\ref{bound}) the square of Pitman slope of $S_{m,n}$ equals the Bahadur local index of
$B_{\emptyset,n}^1$. Thus, the respective left-hand sides in (\ref{ub}) coincide, and
the tests based on
$B_{\emptyset,n}^1$ and $S_{m,n}$, $m\geq 2$, are asymptotically equivalent.}
\end{remark}

\section*{Acknowledgments}
The research of A. Nazarov was partly supported by RFBR grant
10-01-00154a. The research of N. Stepanova was supported by an NSERC
grant. We would like to thank Prof. Yu. V. Tarannikov for
communicating with us regarding references \cite{kl, korsh}. We are
grateful to Prof. Ya. Yu. Nikitin for discussions and suggestions.

\bibliographystyle{gNST}

\end{document}